\documentclass{article}
\begin{document}
\font\germ=eufm10
\def\ssl{\hbox{\germ sl}}
\def\slh{\widehat{\ssl_2}}
\makeatletter
\def\aaa{@}
\centerline{}
\centerline{\Large\bf Irreducible Modules of Finite Dimensional }
\vskip7pt
\centerline{\Large\bf Quantum Algebras of type A at Roots of Unity}
\vskip15pt
\centerline{Toshiki Nakashima}
\vskip10pt
\centerline{Department of Mathematics,}
\centerline{Sophia University, Tokyo 102-8554, JAPAN}
\centerline{e-mail:\,\,toshiki@mm.sophia.ac.jp}
\makeatother

\renewcommand{\labelenumi}{$($\roman{enumi}$)$}
\renewcommand{\labelenumii}{$(${\rm \alph{enumii}}$)$}
\font\germ=eufm10

\def\AA{{\cal A}}
\def\al{\alpha}
\def\beneme{\begin{enumerate}}
\def\beq{\begin{equation}}
\def\beqn{\begin{eqnarray}}
\def\beqnn{\begin{eqnarray*}}
\def\bigsl{{\hbox{\fontD \char'54}}}
\def\bbra#1,#2,#3{\left\{\begin{array}{c}\hspace{-5pt}
#1;#2\\ \hspace{-5pt}#3\end{array}\hspace{-5pt}\right\}}
\def\cd{\cdots}
\def\CC{\hbox{\bf C}}
\def\ddd{\hbox{\germ D}}
\def\del{\delta}
\def\Del{\Delta}
\def\ei{e_i}
\def\eit{\tilde{e}_i}
\def\eneme{\end{enumerate}}
\def\ep{\epsilon}
\def\eeq{\end{equation}}
\def\eeqn{\end{eqnarray}}
\def\eeqnn{\end{eqnarray*}}
\def\fit{\tilde{f}_i}
\def\ft{\tilde{f}}
\def\gau#1,#2{\left[\begin{array}{c}\hspace{-5pt}#1\\
\hspace{-5pt}#2\end{array}\hspace{-5pt}\right]}
\def\ge{\hbox{\germ g}}
\def\gl{\hbox{\germ gl}}
\def\hom{{\hbox{Hom}}}
\def\ify{\infty}
\def\io{\iota}
\def\kp{k^{(+)}}
\def\km{k^{(-)}}
\def\llra{\relbar\joinrel\relbar\joinrel\relbar\joinrel\rightarrow}
\def\lan{\langle}
\def\lar{\longrightarrow}
\def\lm{\lambda}
\def\Lm{\Lambda}
\def\mapright#1{\smash{\mathop{\longrightarrow}\limits^{#1}}}
\def\mm{{\bf{\rm m}}}
\def\nd{\noindent}
\def\nn{\nonumber}
\def\nnn{\hbox{\germ n}}
\def\oint{{\cal O}_{\rm int}(\ge)}
\def\ot{\otimes}
\def\op{\oplus}
\def\opi{\ovl\pi_{\lm}}
\def\ovl{\overline}
\def\plm{\Psi^{(\lm)}_{\io}}
\def\qq{\qquad}
\def\q{\quad}
\def\qed{\hfill\framebox[2mm]{}}
\def\QQ{\hbox{\bf Q}}
\def\qi{q_i}
\def\qii{q_i^{-1}}
\def\ran{\rangle}
\def\rlm{r_{\lm}}
\def\ssl{\hbox{\germ sl}}
\def\slh{\widehat{\ssl_2}}
\def\ti{t_i}
\def\tii{t_i^{-1}}
\def\til{\tilde}
\def\tt{{\hbox{\germ{t}}}}
\def\ttt{\hbox{\germ t}}
\def\ua{U_{\AA}}
\def\ue{U_{\vep}}
\def\uq{U_q(\ge)}
\def\ufin{U^{\rm fin}_{\vep}}
\def\ufinp{(U^{\rm fin}_{\vep})^+}
\def\ufinm{(U^{\rm fin}_{\vep})^-}
\def\ufinz{(U^{\rm fin}_{\vep})^0}
\def\uqm{U^-_q(\ge)}
\def\uqp{U^+_q(\ge)}
\def\uqmq{{U^-_q(\ge)}_{\bf Q}}
\def\uqpm{U^{\pm}_q(\ge)}
\def\uqq{U_{\bf Q}^-(\ge)}
\def\uqz{U^-_{\bf Z}(\ge)}
\def\ures{U^{\rm res}_{\AA}}
\def\urese{U^{\rm res}_{\vep}}
\def\uresez{U^{\rm res}_{\vep,\ZZ}}
\def\util{\widetilde\uq}
\def\vep{\varepsilon}
\def\vp{\varphi}
\def\vpi{\varphi^{-1}}
\def\VV{{\cal V}}
\def\xii{\xi^{(i)}}
\def\Xiioi{\Xi_{\io}^{(i)}}
\def\WW{{\cal W}}
\def\wtil{\widetilde}
\def\what{\widehat}
\def\wpi{\widehat\pi_{\lm}}
\def\ZZ{\hbox{\bf Z}}

\vskip40pt
\renewcommand{\thesection}{\arabic{section}}
\section{Introduction}
\setcounter{equation}{0}
\renewcommand{\theequation}{\thesection.\arabic{equation}}

The quantum group $\uq$ associated with a simple Lie algebra $\ge$
is an associative algebra over the rational function field $\CC(q)$
($q$ is an indeterminate) and we can define its "integral" form over the 
Laurant polynomial ring $\CC[q,q^{-1}]$, which enables us to specialize
$q$ to any non-zero complex number $\vep$. We are going to see 
two types of such integral forms and accordingly, 
we obtain two types of specializations, one is called 
the 'restriced specialization' denoted by $\urese$, 
and the other is called the 'non-restricted specialization'
denoted by $\ue$. Both coincide if $\vep$ is trancendental.
But we are interested in the case that 
$\vep$ is the $l$-th primitive root of unity, 
where $l$ is an odd integer greater than 1. 
In the case, they do not so.
The former is initiated by Lusztig \cite{L1},\cite{L2} and the
latter is introduced in \cite{DK} by DeConcini and Kac.
Their representation theories are quite different:
Irreducible $\urese$-modules are highest weight modules
in some sense and 
the classification of the irreducible modules 
is same as the one for simple Lie algebras or ordinary quantum
algebras
(see Theorem \ref{class} below).  Furthermore, irreducble modules
possess the remarkable property ``tensor product theorem''
(see Theorem \ref{tensor} below), which claims that 
arbitray irreducible highest weight module $V(\lm)$ with the 
highest weight $\lm$ is devided into tensor product of 
two irreducible modules
$V(\lm^{(0)})$ and $V(l\lm^{(1)})$ where $\lm^{(0)}$ and $\lm^{(1)}$
are as in Theorem \ref{tensor}. 
Here the module $V(\lm^{(0)})$ is identified with the 
irreducible $\ufin$-module, wehre $\ufin$ is 
some finite dimenstional subalgebra of $\urese$ (see 2.2)
and the module $V(l\lm^{(1)})$ can be 
identified with the irreducible highest weight $U(\ge)$-module
$V(\lm^{(1)})$, whose structure is known very well. Thus, if the
structure of $V(\lm^{(0)})$ is clarified, we can analize the detailed 
feature of $V(\lm)$. 
Indeed, the character of $V(\lm)$ is given by the famous
Kazhdan-Lusztig formula.
But structures as a module, {\it e.g.}, explicit descriptions of basis vectors
or actions of the generators on them, 
are not still clear.

On the other hand, irreducible $\ue$-modules 
are not necessarily highest or lowest weight
modules.
They are characterized by many continious parameters and 
if they are ``generic'', their  dimensions are all same
(see \cite{DK},\cite{CP}).
But if we specialize the parameters properly, the modules become 
reducible. 
In \cite{DJMM}, Date, Jimbo, Miki and Miwa constrcuted  
such $\ue$-modules for $A_n$-type explicitly, which is called the 
'maximal cyclic representations' that is realized in the vector space
$\VV:=(\CC^l)^{{1\over2}n(n+1)}$.
They contains the continious parameters 
and it is shown that if those parameters are generic, they are irreducible.
Here we consider certain non-generic specialization of the parameters so that 
$\VV$ becomes a reducible $\ue$-module. 
Moreover, we shall observe that 
such a module includes the unique primitive vector (see Proposition
\ref{highest}). 
The submodule generated by this primitive vector can be seen as 
an irreducible $\ufin$-module and isomorphic to $V(\lm^{(0)})$ 
for some $\lm$ (Theorem \ref{main}).

The organizations of the paper is as follows;
in section 2 we review the quantum algebras at roots of unity.
In section 3, we see the maximal cyclic representations 
of the $A$-type following \cite{DJMM} and review
the represenation theory of $\urese$.
In section 4, we specilize the parameters properly and show that 
under the specialization, there exists a unique primitive vector 
in the module.
Finally, in section 5, it is shown that 
the representation space become the module of the
finite dimensional algebra $\ufin$ and 
the submodule generated by the primitive vector is irreducible.

\renewcommand{\thesection}{\arabic{section}}
\section{Algebras at roots of unity}
\setcounter{equation}{0}
\renewcommand{\theequation}{\thesection.\arabic{equation}}

In this section, we review the algebras treated in this article.
\subsection{Restricted integral forms and specializations}

Let $\CC(q)$ be the rational function field in an indeterminate $q$
and denote the ring $\CC[q,q^{-1}]$ by $\AA$.
We use the notations:
$$
[a]_q:={{q^a-q^{-a}}\over{q-q^{-1}}},\q 
[a]_q!:=[a]_q[a-1]_q\cd [2]_q[1]_q, \q
{\gau{m},{k}}_q:={{[m]_q!}\over{[k]_q![m-k]_q!}}.
$$
Let $I:=\{1,2,\cd,n\}$ be the index set and 
$(a_{ij})_{i,j\in I}$ be the Cartan matrix of type A, i.e.,
$a_{ii}=2$ ($1\leq i\leq n$), $a_{i i+1}=a_{i+1 i}=-1$ ($1\leq i\leq n-1$), 
and $a_{ij}=0$ otherwise.
Let us denote the set of roots (resp. positive roots) by $\Delta$ (resp. $\Delta_+$).
Let $\{h_i\}_{i\in I}$ be the set of simple coroots and 
$\{\al_i\}_{i\in I}$ the set of simple roots.
Define the weight lattice $P:=\{\lm\,|\,\lan h_i,\lm\ran\in \ZZ\}$
(resp. the set of dominant integral weights 
$P_+:=\{\lm\,|\,\lan h_i,\lm\ran\in \ZZ_{\geq0}\}$).
Let $\{\Lm_i\}_{i\in I}$ be the fundamental weights which satisfy
$\lan h_i,\Lm_j\ran=\del_{ij}$ and then $P=\oplus_i\ZZ\Lm_i$.
Let $W$ be the Weyl group of type $A_n$, which is generated by the 
simple reflections $s_i$ $(i\in I)$.
The {\it quantum algebra } $\uq$ is the associative algebra generated by  
$e_i, f_i, t^{\pm}_i$ ($i\in I$) and the relations
\begin{eqnarray}
& t_it_i^{-1}=t_i^{-1}t_i=1,\q t_it_j=t_jt_i, \\
& t_ie_jt_i^{-1}=q^{a_{ij}}e_j,\\
& t_if_jt_i^{-1}=q^{-a_{ij}}e_j,\\
& e_if_j-f_je_i={{t_i-t_i^{-1}}\over{q-q^{-1}}},\label{ef-fe}\\
& \hspace{-5mm}\sum_{k=0}^{1-a_{ij}}(-1)e_i^{(k)}e_je_i^{(1-a_{ij}-k)}=
\sum_{k=0}^{1-a_{ij}}(-1)f_i^{(k)}f_jf_i^{(1-a_{ij}-k)}=0\q(i\ne j),
\end{eqnarray}
where 
$e_i^{(k)}:=e_i^k/[k]_q!$ and $f_i^{(k)}:=f_i^k/[k]_q!$.

\noindent
Here we set
$$
{\gau{t_i,p},{r}}_q:=\prod_{s=1}^{r}{{t_iq^{p+1-s}-t_i^{-1}q^{s-p-1}}\over{q^s-q^{-s}}}.
$$
The algebra $\ures$ is the $\AA$-subalgebra of $\uq$ generated by 
$e_i^{(k)}$, $f_i^{(k)}$, $t_i^{\pm}$ and  $\gau{t_i,p},{k}$ 
($i\in I$, $p,k\in \ZZ$ and $k\geq0$), which is called the {\it restricted
integral form}.

Here we can define the {\it restricted specializations} for any
$\vep\in \CC^{\times}$;
\beq
\urese:=\ures\otimes_{\AA}\CC_{\vep},
\eeq
where $\AA$ acts on $\CC_{\vep}:=\CC$ by $f(q)c:=f(\vep)c$ ($c\in \CC$).

\subsection{Finite dimensional quantum algebra}
For $\vep\in \CC^{\times}$ we use the notation
$$
[a]:={{\vep^a-\vep^{-a}}\over{\vep-\vep^{-1}}},\qq
[a]!:=[a][a-1]\cd[2][1],
\gau{m},{k}:=\gau{m},{k}_{q=\vep}.
$$
Since ${\gau{m},{k}}_q\in \CC[q,q^{-1}]$,
the definition of $\gau{m},{k}$ is valid.

As for the specializations of $q$, we shall be  
interested in the case that $\vep$  is a root of unity.
So in what follows, suppose that:
$$
{\hbox{$l$ is the odd integer greater than 1 and 
$\vep$ is the primitive $l$-th root of unity.}}
$$

Under this setting, we can find an interesting finite dimensional 
subalgebra $\ufin$ of $\urese$.
$\ufin$ is defined as the subalgebra of $\urese$ generated by 
$e_i$, $f_i$ and $t_{i}^{\pm}$
($1\leq i\leq n$).
We know that this algebra is finite dimensional over $\CC$ 
with the dimension $2^nl^{n^2+2n}$
(see Proposition \ref{PBW} below).

This $\ufin$ is also defined by "generators and relations" as follows;
\newtheorem{pro2}{Proposition}[section]
\begin{pro2}[\cite{CP},\cite{L1}]
\label{ufin}
The algebra $\ufin$ is isomorphic to the associative $\CC$-algebra
with generators $e_{\al}$, $f_{\al}$ $(\al\in\Delta_+)$ and 
$t_i^{\pm}$ $(1\leq i\leq n)$ satisfying the following relations;
\beqn
& t_it_i^{-1}=t_i^{-1}t_i=1,\q t_it_j=t_jt_i,\label{rel1} \\
& t_ie_jt_i^{-1}=\vep^{a_{ij}}e_j,\label{rel2}\\
& t_if_jt_i^{-1}=\vep^{-a_{ij}}e_j,\label{rel3}\\
& e_if_j-f_je_i=\delta_{ij}{{t_i-t_i^{-1}}\over{\vep-\vep^{-1}}}.\label{rel4}
\eeqn
If $(\al_i,\al)=0$ and $i<g(\al)$,
\beqn
& e_{i}e_{\al}=e_{\al}e_{i}, \label{com-e}\\
& f_{i}f_{\al}=f_{\al}f_{i}.\label{com-f}
\eeqn
If $(\al_i,\al)=-1$ and $i<g(\al)$, 
\beqn
& e_{\al+\al_i}=\vep^{-1} e_{\al}e_i-e_ie_{\al},\label{h-rt+}\\
& \vep e_ie_{\al+\al_i}=e_{\al+\al_i}e_i,\label{com1-e}\\
& \vep e_{\al+\al_i}e_{\al}=e_{\al}e_{\al+\al_i},\label{com2-e}\\
& f_{\al+\al_i}=\vep f_{\al}f_i-f_if_{\al},\label{h-rt-}\\
& \vep f_if_{\al+\al_i}=f_{\al+\al_i}f_i,\label{com1-f}\\
& \vep f_{\al+\al_i}f_{\al}=f_{\al}f_{\al+\al_i}.\label{com2-f}\\
& e_{\al}^l=f_{\al}^l=0, ({\hbox{ for any }}\al\in \Delta_+),\label{00}\\
& t_i^{2l}=1 ({\hbox{ for any }}i\in I),\label{11}
\eeqn
where we define $g(\al)$ $(\al\in\Delta_+)$ to be the largest index
satisfying $c_i\ne0$ if we write $\al=\sum_i c_i\al_i$
and set $e_i:=e_{\al_i}$ and $f_i:=f_{\al_i}$.
\end{pro2}

Define $\ufinp$ (resp. $\ufinm$, $\ufinz$) to be 
the subalgebra of $\ufin$ generated by $e_i$ (resp. $f_i$, $t^{\pm}_i$).
Fix a reduced expression $w_0=s_{i_1}s_{i_2}\cd s_{i_N}$ of the 
longest elememtn of the Weyl group $W$ and 
set $\beta_k:=s_{i_1}s_{i_2}\cd s_{i_{k-1}}(\al_{i_k})$ for $k\in\{1,\cd,N\}$
where $N:={1\over2}n(n+1)$ is the number of positive roots.
Here we have the following Poincar$\acute{\rm e}$-Birkhoff-Witt type theorem:
\begin{pro2}[\cite{CP},\cite{L1}]
\label{PBW}
\begin{enumerate}

\item
The algebra $\ufinp$ is a finite dimensional $\CC$-vector space 
with the basis
\beq
\{e_{\beta_{N}}^{r_N}e_{\beta_{N-1}}^{r_{N-1}}\cd 
e_{\beta_{1}}^{r_1}\}_{0\leq r_1,\cd,r_N<l}.
\eeq
\item
The algebra $\ufinm$ is a finite dimensional $\CC$-vector space 
with the basis
\beq
\{f_{\beta_{N}}^{r_N}f_{\beta_{N-1}}^{r_{N-1}}\cd 
f_{\beta_{1}}^{r_1}\}_{0\leq r_1,\cd,r_N<l}.
\eeq
\item
The algebra $\ufinz$ is a finite dimensional $\CC$-vector space 
with the basis
\beq
\{t_n^{r_n}t_{n-1}^{r_{n-1}}\cd t_{1}^{r_1}\}_{0\leq r_1,\cd,r_n<2l}.
\eeq
\item
Multiplication defines an isomorphism of $\CC$-vector space;
\beq
\ufinm\ot \ufinz\ot \ufinp\mapright{\sim}\ufin.
\eeq
\end{enumerate}
\end{pro2}

\subsection{Non-restrcited specializations}
Here we see another type of specialization of $q$ to a 
root of unity.

Introduce the elements
$$
[t_i;m]:={{t_iq^m-t_i^{-1}q^{-m}}\over{q-q^{-1}}}\in \uq.
$$
The algebra $\ua$ is the $\AA$-subalgebra of $\uq$
generated by the elements $e_i$, $f_i$, $t_i^{\pm}$ and 
$[t_i;0]$ $(1\leq i\leq n)$.

{\sl Remark.\,}
The defining relations for $\ua$ are as in 2.1, but repalcing 
(\ref{ef-fe}) by
\beq
e_if_j-f_je_i=\del_{ij}[t_i;0].
\eeq
and add the relation $(q-q^{-1})[t_i;0]=t_i-t_i^{-1}$.

Now for arbitrary $\vep\in\CC^{\times}$ we define the $\CC$-algebra
$$
\ue:=\ua\ot_{\AA}\CC_{\vep},
$$
where $\AA$ acts on $\CC_{\vep}=\CC$ by $f(q)c=f(\vep)c$ $(c\in\CC)$.
This $\ue$ is called the {\it non-restricted specialization}.

\renewcommand{\thesection}{\arabic{section}}
\section{Representations}
\setcounter{equation}{0}
\renewcommand{\theequation}{\thesection.\arabic{equation}}

\subsection{Maximal cyclic representations of $\ue$}

The representation theory of $\ue$ is discussed in 
\cite{DK} in which the maximal dimension of irreducible representations
for $A_n$ type
is given by $l^{{1\over 2}n(n+1)}$ in the case $\vep$ is the $l$-th
root of unity  and in \cite{DJMM}, 
it is constructed explicitly and called the 'maximal cyclic representaions'. 
Here we modify the presnetations in \cite{DJMM}
subtly in order to simplify the arguments in the section 4.

Let $H$ be the group generated by $\{x_{ij},z_{ij}\}_{1\leq i\leq j\leq n}$ and 
the center $\vep$
with the relations $z_{ij}x_{ij}=\vep x_{ij}z_{ij}$
and all others commute each other, 
and set $\WW:=\CC[H]$ the group ring of $H$.
For $r:=(r_1,\cd,r_n),$ $s:=(s_1,\cd,s_n)\in (\CC^{\times})^n$, 
we define the map $\vp_{r,s}:\ue\longrightarrow \WW$ by (see \cite{DJMM});
\beqn
&& \vp_{r,s}(e_i):=\sum_{k=i}^{n}x_{i\,k}x_{i\, k+1}\cd x_{i\,n}
    \{r_iz_{i\,k}z_{i\, k-1}z_{i-1\, k-1}^{-1}z_{i+1\, k}^{-1}\},
\label{ei}\\
&& \vp_{r,s}(f_i):=\sum_{k=1}^{i} 
    x_{i+1-k\, n+1-k}^{-1}x_{i+2-k\, n+2-k}^{-1}\cd x_{\,in}^{-1}
\label{fi}\\
&&\qq\qq \times\{s_iz_{i+1-k\, n-k}z_{i+1-k\, n+1-k}^{-1} z_{i-k\, n+1-k} 
z_{i-k\, n-k}^{-1}\}, \nn\\
&& \vp_{r,s}(t_i):={{r_i}\over{s_i}}z_{i\,n}^2z_{i-1\, n}^{-1}z_{i+1\, n}^{-1},
\eeqn
where we use the notation $\{z\}=(z-z^{-1})/({\vep-\vep^{-1}})$.

Let $*:\WW\longrightarrow \WW$ be the $\CC$-linear involution defined by 
$$
x_{jk}^*:=x_{k+1-j \,k}^{-1},\qq
z_{jk}^*:=z_{k+1-j \,k}^{-1}, 
$$
and set
$$
A_{ik}:=x_{i\,k}x_{i\, k+1}\cd x_{i\,n},\qq
B_{ik}:=z_{i\,k}z_{i\, k-1}z_{i-1\, k-1}^{-1}z_{i+1\, k}^{-1}.
$$
Then, (\ref{ei}) and (\ref{fi}) can  be written in the following forms;
\beq
\vp_{r,s}(e_i)=\sum_{k=i}^{n}A_{ik}\{r_i B_{ik}\},
\qq \vp_{r,s}(f_i)=\sum_{k=1}^{i}A_{n+1-i\,n+1-k}^*\{s_i B_{n+1-i\, n+1-k}^*\}.
\eeq

\newtheorem{pro3}{Proposition}[section]
\begin{pro3}
\label{alghom}
The map $\vp_{r,s}$ defines a $\CC$-linear algebra homomorphism
from $\ue$ to $\WW$.
\end{pro3}

\newtheorem{lm3}[pro3]{Lemma}
\begin{lm3}
\label{AB}
The following commutation relations hold (see \cite{DJMM}(2.5)).
\beqnn
A_{ij}B_{ik}&& =\vep^{-2}B_{ik}A_{ij}\qq{\rm if }\q j<k,\\
            && =\vep^{-1}B_{ik}A_{ij}\qq{\rm if }\q j=k,\\
            && =B_{ik}A_{ij}\qq{\rm if }\q j>k.
\eeqnn
\end{lm3}
{\sl Proof of Proposition \ref{alghom}.}
We have $A_{ik}B_{ik}=\vep^{-1}B_{ik}A_{ik}$ and then
$$
\vp_{r,s}(e_i)=\sum_{k=i}^{n}\{r_i\vep^{-1} B_{ik}\}A_{ik},
\qq \vp_{r,s}(f_i)=\sum_{k=1}^{i}\{s_i\vep^{-1} B_{n+1-i\,n+1-k}^*\}A_{n+1-i\,n+1-k}^*.
$$
This implies $\vp_{r,s}=\rho_{\vep^{-1}r,\vep^{-1}s}$ 
($\rho_{r,s}$ is given in \cite{DJMM}).
Thus, by Theorem 2.2 in \cite{DJMM}, we obtained the desired result.\qed

\begin{pro3}
\label{power}
For any $m\in \ZZ_{>0}$, we have
\beqn
&&
\hspace{-20pt}\vp_{r,s}(e_i^m)=[m]!\sum_{p=1}^m 
\sum_{\tiny\begin{array}{c}
 i\leq k_p<\cd < k_1\leq n\\
 1\leq \nu_p<\cd < \nu_{1}=m
\end{array}}
\prod_{r=1}^{p} A_{i \, k_r}^{\nu_r-\nu_{r+1}} 
\prod_{r=1}^{p} \bbra{r_iB_{i \, k_r}},{\nu_{r}-1},{\nu_r-\nu_{r+1}},
\label{em}\\
&& \hspace{-20pt}\vp_{r,s}(f_i^m) \nn\\
&&\hspace{-20pt} =[m]!\sum_{p=1}^m 
\sum_{\tiny\begin{array}{c}
 i\leq k_p<\cd < k_1\leq n\\
 1\leq \nu_p<\cd < \nu_{1}=m
\end{array}}
\prod_{r=1}^p A_{n+1-i \, n+1-k_r}^{*\nu_r-\nu_{r+1}} 
\prod_{r=1}^p \bbra{s_iB_{n+1-i \,n+1-k_r}^*},{\nu_{r}-1},{\nu_r-\nu_{r+1}},\qq
\label{fm}
\eeqn
\end{pro3}
where $\nu_{p+1}=0$ and we set
$$
\bbra{a},{b},{c}:={{\{a\vep^b\}\{a\vep^{b-1}\}\cd\{a\vep^{b-c+1}\}}\over{[c]!}}.
$$
{\sl Remark.\,} 
The definition of $\bbra{a},{b},{c}$ is invalid for
$\vep$ such that $[c]!=0$. But in the right hand-side of 
(\ref{em}) and (\ref{fm})  we see that the term
$$
{{[m]!}\over{\prod_{r=1}^p [\nu_r-\nu_{r+1}]}}\qq
(1\leq \nu_p<\cd < \nu_{1}=m).
$$
is valid since ${{[m]_q!}/{\prod_{r=1}^p [\nu_r-\nu_{r+1}]_q}}\in\ZZ[q,q^{-1}]$.

{\sl Proof.\,}
In \cite{DJMM}, the following formula is given
\beq
\rho_{r,s}(e_i^m)=[m]!\sum_{p=1}^m 
\sum_{\tiny\begin{array}{c}
 i\leq k_p<\cd < k_1\leq n\\
 1\leq \nu_p<\cd < \nu_{1}=m
\end{array}}
\prod_{r=1}^{p} \bbra{r_iB_{i \, k_r}},{-\nu_{r+1}},{\nu_r-\nu_{r+1}}
\prod_{r=1}^{p} A_{i \, k_r}^{\nu_r-\nu_{r+1}},\label{djmm-em}
\eeq
where $\nu_{p+1}=0$.
Since $\vp_{r,s}=\rho_{\vep^{-1}r,\vep^{-1}s}$, 
it follows from (\ref{djmm-em}) 
\beqn
&& \hspace{-15pt}\vp_{r,s}(e_i^m)=\rho_{\vep^{-1}r,\vep^{-1}s}(e_i^m)\nn\\
&& \hspace{-15pt} =[m]!\sum_{p=1}^m 
\sum_{\tiny\begin{array}{c}
 i\leq k_p<\cd < k_1\leq n\\
 1\leq \nu_p<\cd < \nu_{1}=m
\end{array}}
\prod_{r=1}^{p} \bbra{\vep^{-1}r_iB_{i \, k_r}},{-\nu_{r+1}},{\nu_r-\nu_{r+1}}
\prod_{r=1}^{p} A_{i \, k_r}^{\nu_r-\nu_{r+1}}
\label{em2}
\eeqn

Here by Lemma \ref{AB} we obtain for $ i\leq k_p<\cd < k_1\leq n$ and
$1\leq r\leq p$
\beqnn
\hspace{-20pt}\bbra{\vep^{-1}r_iB_{i \, k_r}},{-\nu_{r+1}},{\nu_r-\nu_{r+1}}
\left(\prod_{r=1}^{p} A_{i \, k_r}^{\nu_r-\nu_{r+1}}\right)
& = &
\left(\prod_{r=1}^{p} A_{i \, k_r}^{\nu_r-\nu_{r+1}}\right)
\bbra{\vep^{-1+\nu_{r}+\nu_{r+1}-2\nu_{p+1}}r_iB_{i \, k_r}},{-\nu_{r+1}},{\nu_r-\nu_{r+1}}\\
& = &
\left(\prod_{r=1}^{p} A_{i \, k_r}^{\nu_r-\nu_{r+1}}\right)
\bbra{r_iB_{i \, k_r}},{\nu_{r}-1},{\nu_r-\nu_{r+1}},
\eeqnn
where we use $\nu_{p+1}=0$.
Thus, we obtain (\ref{em}). Similarly we also get (\ref{fm}).\qed

\vskip10pt

Let $V_{i\,j}$ $(1\leq i\leq j\leq n)$ be a copy of the vector space $\CC^l$ 
and set $\VV:=\otimes_{1\leq i\leq j\leq n}V_{i\,j}$.
Let $u_0,\cd,u_{l-1}$ be the standard basis of $\CC^l$.
Now we define the representation $(\psi_{a,b},\VV)$ of $\WW$ as follows:
Let  $Z_{j\,k},X_{j\,k}\in {\rm End}(\VV)$ be the matrices 
defined as $Z_{j\,k}u_i=u_{i+1}$ and $X_{j\,k}u_{i}=\vep^iu_i$ on the 
component $V_{jk}$ and as the identity on the other component.
For non-zero parameters $a:=(a_{i\,j})_{1\leq i\leq j\leq n}$ and 
$b:=(b_{i\,j})_{1\leq i\leq j\leq n}\in (\CC^{\times})^{n(n+1)/2}$, define 
$\psi_{a,b}(x_{i\,j}), \,\psi_{a,b}(z_{i\,j})\in{\rm End}(\VV)$
to be 
\beq
\psi_{a,b}(x_{ij})=a_{ij}X_{ij},\qq
\psi_{a,b}(z_{ij})=b_{ij}Z_{ij}.
\eeq
We can easily check that these define the representation of $\WW$:
\beq
\psi_{a,b}:\WW\longrightarrow {\rm End}(\VV).
\eeq
Composing $\vp_{r,s}$ and $\psi_{a,b}$;
$$
\Phi_{r,s,a,b}:=\psi_{a,b}\circ\vp_{r,s}
\,:\, \ue\mapright{\vp_{r,s}}\WW\mapright{\psi_{a,b}}{\rm End}(\VV).
$$
we obtain the representation of $\ue$ denoted by $(\Phi_{r,s,a,b},\VV)$.
The representation introduced in \cite{DJMM} is just as 
$(\Phi_{\vep r,\vep s,a,b},\VV)$ in our notation
since we have $\vp_{r,s}=\rho_{\vep^{-1}r,\vep^{-1}s}$ in the 
proof of Proposition \ref{alghom}.

In \cite{DJMM}, it is shown that 
the central elements of $\ue$ take
values in an open set of $\CC^{n(n+2)}$ 
and then by \cite{DK}, it turns out to be generically irreducible
for the parameters $r,s,a,b$.
We are interested in specializations of these parameters so that 
the representation $(\Phi_{r,s,a,b}),\VV)$ is not necessarily irreducible.

\subsection{Representations of $\urese$ and $\ufin$}

We review the representation theory of $\urese$.
The classification of the irreducible represntations of $\urese$
is given by Lusztig \cite{L2}. Before seeing it, let us recall the
notions of highest weight modules.
\newtheorem{def3}[pro3]{Definition}

\begin{def3}
Let $V$ be a $\urese$-module
of type ${\bf 1}$ (as for ``type'', see \cite{CP},\cite{L2}).
\begin{enumerate}
\item
The weight spaces $V_{\lm}$ $(\lm=\sum_im_i\Lm_i\in P)$ 
of $V$  are
defined by 
\beq
V_{\lm}:=\left\{v\in V| 
t_iv=\vep^{m_i^{(0)}}v,\q
\gau{t_i;0},{l}v=\gau{m^{(1)}_i},{l}v \right\},
\eeq
where $m_i=m_i^{(0)}+lm_i^{(1)}$ and $0\leq m_i^{(0)}<l$.
\item
$V$ is a highest weight module if $V$ 
is generated by a primitive vector,
i.e., a vector $v\in V_{\lm}$ 
for some $\lm\in P$, such that $e_iv=e^{(l)}_iv=0$ for any $i\in I$.
In the case, $\lm$ is called the highest weight and $v$ 
is called the highest weight vector of $V$.
\end{enumerate}
\end{def3}

Let $V(\lm)$ be the irreducible highest weight $\uq$-module 
given by $V(\lm)=\uq/{\cal I}$ 
where $\lm\in P_+$ and, 
${\cal I}$ is the left ideal generated by
$e_i$, $f_i^{1+\lan h_i,\lm\ran}$ and $t_i-q^{\lan h_i,\lm\ran}$
$(1\leq i\leq n)$. Here denote the generator of $V(\lm)$
by $v_{\lm}$.
Let $V^{\rm res}_{\AA}(\lm)$ be the $\ures$-submodule of $V(\lm)$ 
generated by $v_{\lm}$. 
Set $W^{\rm res}_{\vep}(\lm):=V^{\rm res}_{\AA}(\lm)\ot_{\AA}\CC_{\vep}$,  
which is naturally 
$\urese$-module.
 Note that $W^{\rm res}_{\vep}(\lm)$ is not necessarily irreducible.
So, let $Y$ be its maximal proper submodule and define 
$V^{\rm res}_{\vep}(\lm):=W^{\rm res}_{\vep}(\lm)/Y$ to be the
irreducibel quotient, which is type ${\bf 1}$ highest weight module
with the highest weight $\lm$.

\newtheorem{thm3}[pro3]{Theorem}
\begin{thm3}[\cite{L2}]
\label{class}
Arbitrary finite-dimensional irreducible $\urese$-module $V$ of type {\bf 1}
is isomorphic to $V^{\rm res}_{\vep}(\lm)$ for a unique $\lm\in P_+$.
\end{thm3}
Note that 
arbitrary finite-dimensional irreducible $\urese$-module $V$ of type {\bf 1}
is a direct sum of its weight spaces.

\begin{thm3}[\cite{L2}]
\label{tensor}
For $\lm=\sum_im_i\Lm_i\in P_+$, define $\lm^{(0)}:=\sum_im_i^{(0)}\lm_i$
and $\lm^{(1)}:=\sum_im_i^{(1)}\lm_i$ where $m_i=m_i^{(0)}+lm_i^{(1)}$
with $0\leq m_i^{(0)}<l$ (and then $\lm=\lm^{(0)}+l\lm^{(1)}$).
The $\urese$-module $V^{\rm res}_{\vep}(\lm)$ is 
isomorphic to $V(\lm^{(0)})\ot V(l\lm^{(1)})$.
\end{thm3}

Here we call a weight $\lm\in P_+$ satisfying 
$\lm=\lm^{(0)}$ a $l$-{\it restricted weight}.
As we have stated in the introduction, the module $V(\lm^{(0)})$ is irreducible 
$\ufin$-module and $V(l\lm^{(1)})$ is identified with the irreducible 
highest weight $\ssl_{n+1}$-module $V(\lm^{(1)})$. Since we know the structure 
of irreducible $\ssl_{n+1}$-module well,  this theorem implies that 
the structure of the module $V^{\rm res}_{\vep}(\lm)$ can be clarified
if we shall make clear the one for $V(\lm^{(0)})$.

\renewcommand{\thesection}{\arabic{section}}
\section{Primitive vectors}
\setcounter{equation}{0}
\renewcommand{\theequation}{\thesection.\arabic{equation}}

Let $l$ and $\vep$ be same as in the previous section.
\subsection{Specializations of parameters}

Let $M:=\{\mm=(m_{jk})_{1\leq j\leq k\leq n}|0\leq m_{jk}\leq l-1\}$
be the index set of the standard basis of $\VV$.
We can consider the additive structure on $M$ via
the natural identification $M\cong (\ZZ/l\ZZ)^{{1\over2}n(n+1)}$.
For $\mm\in M$ we write $u_{\mm}:=\otimes_{1\leq j\leq k\leq n}u_{m_{jk}}$
($u_{m_{jk}}\in V_{jk}$).

Here we consider the following specialization of parameters 
${r,s,a,b}$ :
\beqn
&& a_{i\,k}a_{i\,k+1}\cd a_{i\,n}=1,\label{para1}\\
&& r_ib_{i\,k}b_{i\, k-1}b_{i-1\,k-1}^{-1}b_{i+1\, k}^{-1}=1,
(1\leq i\leq k\leq n),\label{para2}\\
&& {{r_i}\over {s_i}}b_{i\, n}^2b_{i-1\, n}^{-1}b_{i+1\, n}^{-1}=\vep^{\lm_i}\q
\label{para3}
\eeqn 
where integers $\{\lm_i\}_{1\leq i\leq n}$ satisfy  $0\leq \lm_i<l$.

{\sl Remark.}\,
Here note that 
the set of parameters satisfying (\ref{para1})--(\ref{para3}) is never empty.
Indeed, if we set $a_{j\,k}=b_{j\,k}=1$ for any $(j,k)$ and 
$r_i=1$ and $s_i=\vep^{-\lm_i}$ for any $i$, it is trivial to see that
these satisfy (\ref{para1})--(\ref{para3}). 
(By (\ref{para1}), we have $a_{jk}=1$ for all $1\leq j\leq k\leq n$).

\newtheorem{lm4}{Lemma}[section]
\begin{lm4}
\label{si}
Under the specialization (\ref{para2}) and
(\ref{para3}), we have
\beq
s_ib_{i+1-k n-k}b_{i+1-k n+1-k}^{-1}b_{i-k n+1-k}b_{i-kn-k}^{-1}=\vep^{-\lm_i}.
\label{bb1}
\eeq
\end{lm4}

{\sl Proof.}
Using (\ref{para2}), we have $\,r_ib_{ik}b_{ik-1}b_{i-1\,k-1}^{-1}b_{i+1\,k}^{-1}=1=
r_ib_{ik-1}b_{ik-2}b_{i-1\,k-2}^{-1}b_{i+1\,k-1}^{-1}$ and then
$b_{ik}b_{i-1\,k-1}^{-1}b_{i+1k}^{-1}=b_{ik-2}b_{i-1\,k-2}^{-1}b_{i+1\,k-1}^{-1}$.
Changing $i\rightarrow i-k$ and $k\rightarrow n-k$, we get
\beq
b_{i+1-k\,n-k}b_{i+1-k\,n+1-k}^{-1}
=b_{i-k\,n-k-1}b_{i-k\,n+1-k}^{-1}b_{i-k-1\,n-k-1}b_{i-k-1\,n-k-1}^{-1}.
\label{bb2}
\eeq
By (\ref{para2}) with $k=n$ and (\ref{para3}), we have 
$s_ib_{i n-1}b_{i n}^{-1}b_{i-1\, n}b_{i-1\,n-1}^{-1}=\vep^{-\lm_i}$, 
which is (\ref{bb1}) in the case $k=1$.
Suppose that (\ref{bb1}) holds and substitute (\ref{bb2}) into (\ref{bb1}).
Then we obtain 
$$
s_ib_{i-k n-k-1}b_{i-k n-k}^{-1}b_{i-k-1 n-k}b_{i-k-1n-k-1}^{-1}=\vep^{-\lm_i}.
$$
Thus, the induction on $k$ procedes
and then we prove (\ref{bb1}) for any $k\in \{1,2,\cd,i\}$.\qed

\vskip10pt
By  (\ref{para1})--(\ref{para3}) and this lemma, we have 
\beqn
&& \Phi_{r,s,a,b}(e_i):=\sum_{k=i}^{n}X_{i\,k}X_{i\, k+1}\cd X_{i\,n}
    \{Z_{i\,k}Z_{i\, k-1}Z_{i-1\, k-1}^{-1}Z_{i+1\, k}^{-1}\},
\label{ei-spe}
\\
&& \Phi_{r,s,a,b}(f_i):=\sum_{k=1}^{i} 
    X_{i+1-k\, n+1-k}^{-1}X_{i+2-k\, n+2-k}^{-1}\cd X_{\,in}^{-1}
\label{fi-spe}
\\
&&\qq\qq \times\{\vep^{-\lm_i}Z_{i+1-k\, n-k}Z_{i+1-k\, n+1-k}^{-1} Z_{i-k\, n+1-k} 
Z_{i-k\, n-k}^{-1}\}, \nn\\
&& \Phi_{r,s,a,b}(t_i):=\vep^{\lm_i}Z_{i\,n}^2Z_{i-1\, n}^{-1}Z_{i+1\, n}^{-1},
\eeqn

\subsection{Primitive vecotrs in $\VV$}
Under the specialization in 4.1, 
we get the following:

\newtheorem{pro4}[lm4]{Proposition}
\begin{pro4}
\label{highest-vec}
Under the specialization (\ref{para1}) and (\ref{para2}),
$v\in\VV$ satisfies the condition  
\beq
e_i v=0 \q {\rm for\,\, any\,\,} i=1,\cd, n,
\label{highest}
\eeq
 if and only if $v=cu_{\vec 0}$ $(c\in \CC)$
where ${\vec 0}=(0,0,\cd,0)\in M$.
\end{pro4}

{\sl Proof.\,}
By (\ref{para1}) and (\ref{para2}), the action of $e_i$ on $u_{\mm}\in \VV$
($\mm=(m_{g\,h})\in M$) is given by
\beq
e_iu_{\mm}=\sum_{i\leq k\leq n}
[m_{i\,k}+m_{i\,k-1}-m_{i-1\,k-1}-m_{i+1\,k}]u_{\mm+\ep_{i\,k}+\cd\ep_{i\,n}},
\label{act-e}
\eeq
where $\ep_{j\,k}\in M$ satisfies that
the $(j,k)$-entry is 1 and all others are 0.
If $\mm=\vec 0$, we have $m_{i\,k}+m_{i\,k-1}-m_{i-1\,k-1}-m_{i+1\,k}=0$
for all $i\leq k\leq n$, which implies that $e_iu_{\vec 0}=0$
for any $i$.

Conversely, assume that $v=\sum_{\mm\in M}c_{\mm}u_{\mm}$ $(c_{\mm}\in\CC)$
satisfies (\ref{highest}).
First, we have
\beq
0=e_nv=X_{n\,n}\{Z_{n-1\,n-1}^{-1}Z_{n\,n}\}v
=\sum_{\mm\in M}c_{\mm}[m_{n\,n}-m_{n-1\,n-1}]u_{\mm+\ep_{n\,n}}.
\eeq
This implies that 
\beq
m_{n-1\,n-1}\ne m_{n\,n}\Longrightarrow c_{\mm}=0,
\label{n}
\eeq
and then we have $v=\sum_{\mm\in M,\,m_{n-1\,n-1}=m_{nn}}c_{\mm}v_{\mm}$.

Next, by $e_{n-1}v=0$ we have
\beqnn
&&\hspace{-15pt}0=e_{n-1}v=(X_{n-1\,n-1}X_{n-1\,n}\{Z_{n-1\,n-1}Z_{n-2\,n-2}^{-1}\}
+X_{n-1\,n}\{Z_{n-1\,n}Z_{n-1\,n-1}Z_{n-2\,n-1}^{-1}Z_{n\,n}^{-1}\})v\\
&&=\sum_{\tiny\begin{array}{c}\mm\in M,\\m_{n-1\,n-1}=m_{n\,n}\end{array}}
c_{\mm}[m_{n-1\,n-1}-m_{n-2\,n-2}]u_{\mm+\ep_{n-1\,n-1}+\ep_{n-1\,n}}\\
&&\qq\qq\qq+c_{\mm}[m_{n-1\,n}+m_{n-1\,n-1}-m_{n-2\,n-1}-m_{n\,n}]u_{\mm+\ep_{n-1\,n}}.
\eeqnn
This implies that 
\beq
c_{\mm}[m_{n-1\,n-1}-m_{n-2\,n-2}]=
c_{\mm}[m_{n-1\,n}+m_{n-1\,n-1}-m_{n-2\,n-1}-m_{n\,n}]=0\label{n-1}
\eeq
for any $\mm\in M$ satisfying $m_{n-1\,n-1}=m_{n\,n}$
since all vectors appear in the summation are linearly independent
under the condition $m_{n-1\,n-1}=m_{n\,n}$, that is, 
the index ${\mm+\ep_{n-1\,n-1}+\ep_{n-1\,n}}$ and 
 ${\mm'+\ep_{n-1\,n}}$ never coincide for arbitrary
$\mm,\mm'$ under the condition $m_{n-1\,n-1}=m_{n\,n}$.
Thus, by  (\ref{n}) and (\ref{n-1}) we have unless 
\beqnn
&& m_{n-2\,n-2}=m_{n-1\,n-1}=m_{n\,n},\\
&& m_{n-2\,n-1}=m_{n-1\,n},
\eeqnn
$c_{\mm}=0$.

Here  we assume that $c_{\mm}=0$ in $v$ unless
\beqn
&& m_{i\,i}=m_{i+1\,i+1}=\cd\cd\cd=m_{n\,n},\nn\\
&& m_{i\,i+1}=m_{i+1\,i+2}=\cd\cd=m_{n-1\,n},\label{i}\\
&&\qq\qq \cd\cd \nn\\
&& m_{i\,n-1}=m_{i+1\,n}.\qq\nn
\eeqn

By $e_iv=0$ we get 
\beqn
&&\hspace{-18pt}
0=\sum_{k=i}^nX_{i\,k}X_{i\,k+1}\cd X_{i\,n}
\{Z_{i\,k}Z_{i\,k-1}Z_{i-1\,k-1}^{-1}Z_{i+1\,k}^{-1}\}v\nn \\
&& =\sum_{\tiny\begin{array}{c}\mm\in M,\\ 
\mm{\hbox{ satisfies (\ref{i})}}\end{array}}\sum_{k=i}^n 
c_{\mm}[m_{i\,k}+m_{i\,k-1}-m_{i-1\,k-1}-m_{i+1\,k}]
u_{\mm+\ep_{i\,k}+\cd+\ep_{i\,n}}\nn\\
&& =\sum_{\tiny\begin{array}{c}\mm\in M,\nn\\ 
\mm{\hbox{ satisfies (\ref{i})}}\end{array}}
(c_{\mm}[m_{i\,i}-m_{i-1\,i-1}]u_{\mm+\ep_{i\,i}+\cd+\ep_{i\,n}}\nn\\
&&\qq\qq +c_{\mm}[m_{i\,i+1}+m_{i\,i}-m_{i-1\,i}-m_{i+1\,i+1}]
   u_{\mm+\ep_{i\,i+1}+\cd+\ep_{i\,n}} \label{sum-i}\\
&&\qq\qq\qq\qq\qq\qq +\cd\cd+\nn\\
&&\qq\qq
+c_{\mm}[m_{i\,n}+m_{i\,n-1}-m_{i-1\,n-1}-m_{i+1\,n}]u_{\mm+\ep_{i\,n}}).\nn
\eeqn
It follows from (\ref{i}) that 
all vectors apear in the summation (\ref{sum-i}) are linearly
independent. Therefore, we obtain  that $c_{\mm}=0$,
unless
\beqnn
&& m_{i\,i}-m_{i-1\,i-1}=0,\\
&& m_{i\,i+1}+m_{i\,i}-m_{i-1\,i}-m_{i+1\,i+1}=0,\\
&& \cd\cd\cd\\
&& m_{i\,n}+m_{i\,n-1}-m_{i-1\,n-1}-m_{i+1\,n}=0.
\eeqnn
Thus, from this and (\ref{i}) we get $c_{\mm}=0$ unless
\beqn
&& m_{i-1\,i-1}=m_{i\,i}=m_{i+1\,i+1}=\cd\cd\cd=m_{n\,n},\nn\\
&& m_{i-1\,i}=m_{i\,i+1}=m_{i+1\,i+2}=\cd\cd=m_{n-1\,n},\label{i-1}\\
&&\qq\qq \cd\cd \nn\\
&& m_{i-1,n-2}=m_{i\,n-1}=m_{i+1\,n},\nn\\
&& m_{i-1\,n-1}=m_{i\,n}.\nn
\eeqn

Thus, using $e_2v=e_3v=\cd=e_{n-1}v=e_nv=0$ we have $c_{\mm}=0$
unless
\beqn
&& m_{1\,1}=m_{2\,2}=m_{3\,3}=\cd\cd\cd=m_{n\,n},\nn\\
&& m_{1\,2}=m_{2\,3}=m_{3\,4}=\cd=m_{n-1\,n},\label{2}\\
&&\qq\qq \cd\cd \nn\\
&& m_{1,n-2}=m_{2\,n-1}=m_{3\,n},\nn\\
&& m_{1\,n-1}=m_{2\,n}.\nn
\eeqn
Finally, using $e_1v=0$, we have
\beqn
&&\hspace{-18pt}
0=\sum_{k=1}^nX_{1\,k}X_{1\,k+1}\cd X_{1\,n}
\{Z_{1\,k}Z_{1\,k-1}Z_{2\,k}^{-1}\}v\nn \\
&& =\sum_{\tiny\begin{array}{c}\mm\in M,\\ 
\mm{\hbox{ satisfies (\ref{2})}}\end{array}}\sum_{k=1}^n 
c_{\mm}[m_{1\,k}+m_{1\,k-1}-m_{2\,k}]
u_{\mm+\ep_{1\,k}+\cd+\ep_{1\,n}}\nn\\
&& =\sum_{\tiny\begin{array}{c}\mm\in M,\nn\\ 
\mm{\hbox{ satisfies (\ref{2})}}\end{array}}
(c_{\mm}[m_{1\,1}]u_{\mm+\ep_{1\,1}+\cd+\ep_{1\,n}}\nn\\
&&\qq\qq +c_{\mm}[m_{1\,2}+m_{1\,1}-m_{2\,2}]
   u_{\mm+\ep_{1\,2}+\cd+\ep_{1\,n}} \label{sum-1}\\
&&\qq\qq\qq\qq\qq\qq +\cd\cd+\nn\\
&&\qq\qq
+c_{\mm}[m_{1\,n}+m_{1\,n-1}-m_{2\,n}]u_{\mm+\ep_{1\,n}}).\nn
\eeqn
Under the condition of (\ref{2}), we get that 
$c_{\mm}=0$, unless
\beq
0=m_{1\,1}=m_{1\,2}=\cd=m_{1\,n}.\label{1}
\eeq
Therefore, it follows from (\ref{2}) and (\ref{1}) that $c_{\mm}=0$ unless
\beqn
&& 0=m_{1\,1}=m_{2\,2}=m_{3\,3}=\cd\cd\cd=m_{n\,n},\nn\\
&& 0=m_{1\,2}=m_{2\,3}=m_{3\,4}=\cd=m_{n-1\,n},\label{0}\\
&&\qq\qq \cd\cd \nn\\
&& 0=m_{1,n-2}=m_{2\,n-1}=m_{3\,n},\nn\\
&& 0=m_{1\,n-1}=m_{2\,n},\nn
\eeqn
which implies that $v=cu_{\vec 0}$.\qed

\vskip10pt
{\sl Remark.\,}
In the proof of the proposition, we see easily that we do not need
(\ref{para1}) essentially.
It is required for simplification of the proof or the presentations.
So it is possible to proceed the same argument for generic $a_{jk}$'s.

\vskip 10pt 

The primitive vector $u_{\vec 0}$ possesses the following property:
\begin{pro4}
\label{fi-lm}
Under the condition (\ref{para1}), (\ref{para2}) and (\ref{para3}), we have
$f_i^{\lm_i+1}u_{\vec 0}=0$.
\end{pro4}

{\sl Proof.}
We obtain the explicit form of $f_i^{\lm_i+1}$ on $\VV$ by
(\ref{fm}) in Proposition \ref{power} taking 
$m=\lm_i+1$.
By Lemma \ref{si}, 
under the specialization (\ref{para1}), (\ref{para2}) and
(\ref{para3}), we have
$s_ib_{i+1-k n-k}b_{i+1-k n+1-k}^{-1}b_{i-k n+1-k}b_{i-kn-k}^{-1}=\vep^{-\lm_i}$.
Thus, on $u_{\vec 0}$ we have
\beqnn
&&f_i^{\lm_i+1}u_{\vec 0}\\
&&=
[\lm_i+1]!\sum_{p=1}^{\lm_i+1} 
\sum_{\tiny\begin{array}{c}
 i\leq k_p<\cd < k_1\leq n\\
 1\leq \nu_p<\cd < \nu_{1}=\lm_i+1
\end{array}}
\prod_{r=1}^p A_{n+1-i \, n+1-k_r}^{*\nu_r-\nu_{r+1}} 
\prod_{r=1}^p \bbra{\vep^{-\lm_i}},{\nu_{r}-1},{\nu_r-\nu_{r+1}}u_{\vec 0}.
\eeqnn
For any $p\in\{1,2,\cd,\lm_i+1\}$ and any 
$\nu_p,\cd,\nu_1$ satisying $1\leq \nu_p<\cd<\nu_1=\lm_i+1$, there exists some
$r\in\{1,2,\cd,p\}$
such that $\nu_{r+1}\leq \lm_i<\nu_{r}$.
Therefore, for such $r$ we have
$$
\bbra{\vep^{-\lm_i}},{\nu_{r}-1},{\nu_r-\nu_{r+1}}
={{\{\vep^{-\lm_i+\nu_{r}-1}\}\{\vep^{-\lm_i+\nu_{r}-2}\}\cd
\{\vep^{-\lm_i+\nu_{r+1}}\}}\over{[\nu_r-\nu_{r+1}]!}}
=0.
$$
Thus, we obtain $f_i^{\lm_i+1}u_{\vec 0}=0$.\qed

\vskip10pt

Here for $\lm:=(\lm_1,\cd,\lm_n)$ $(0\leq \lm_i<l)$
we define the $\ue$-submodule $L(\lm)$ of $\VV$ 
by $L(\lm):=\ue u_{\vec 0}$.

Let $V(\lm)$ be the irreducible highest weight $\uq$-module 
as in 3.2.
Let $V_{\AA}(\lm)$ be the $\ua$-submodule of $V(\lm)$ 
generated by $v_{\lm}$. 
Set $V_{\vep}(\lm):=V_{\AA}(\lm)\ot_{\AA}\CC_{\vep}$,  which is naturally 
$\ue$-module.
 Note that $V_{\vep}(\lm)$ is not necessarily irreducible.
By Propositoin \ref{fi-lm}, $f_i^{\lan h_i,\lm\ran+1}u_{\vec 0}=0$ ($i\in I$), thus
we have the following surjective $\ue$-linear map
$\pi: V_{\vep}(\lm) \longrightarrow  L(\lm)$ given by 
$\pi:v_{\lm}\mapsto u_{\vec 0}$.
It seems that the module $L(\lm)$ is in the similar stream of 
the theory of $\urese$-modules.
Here we expect that $L(\lm)$ is an irreducible highest weight $\ue$-module.
Suprisingly, in the next section we obtain more interesting results that 
$L(\lm)$ can be seen as an irreducible $\ufin$-module.
This means that  we get  $\urese$
(or $\ufin$)-modules directly from $\ue$-modules.

\subsection{Shifts of parameters}
Let $r^{(0)}=(r^{(0)}_{j})$, $s^{(0)}=(s^{(0)}_{j})\in (\CC^{\times})^n$ and
$a^{(0)}=(a^{(0)}_{jk})$, $b^{(0)}=(b^{(0)}_{jk})\in(\CC^{\times})^N$ 
be the parameters satisfying (\ref{para1})--(\ref{para3}).

Fix a basis vector $u_{\xi}\in\VV$ ($\xi=(\xi_{jk})\in M$) arbitrarily
and set $b^{(\xi)}:=(\vep^{-\xi_{jk}}b^{(0)}_{jk})\in(\CC^{\times})^N$.
In this setting we obtain the following:
\begin{pro4}
\label{shift}
For any $\mu=(\mu_{jk})\in M$ and any $X\in \ue$,
set 
\beq
\Phi_{r^{(0)},s^{(0)},a^{(0)},b^{(0)}}(X)u_{\mu}=\sum_{\mm\in M}C_{\mm}u_{\mm}.
\eeq
Then we have
\beq
\Phi_{r^{(0)},s^{(0)},a^{(0)},b^{(\xi)}}(X)u_{\mu+\xi}
=\sum_{\mm\in M}C_{\mm}u_{\mm+\xi}
\eeq
\end{pro4}

{\sl Proof.\,}
It is shown easily from the formula
\beqnn
&&\psi_{a^{(0)},b^{(0)}}(x_{jk})u_{\mu_{jk}}=a^{(0)}_{jk}u_{\mu_{jk}+1},\qq
  \psi_{a^{(0)},b^{(0)}}(z_{jk})u_{\mu_{jk}}=b^{(0)}_{jk}\vep^{\mu_{jk}}u_{\mu_{jk}},\\
&&\psi_{a^{(0)},b^{(\xi)}}(x_{jk})u_{\mu_{jk}+\xi_{jk}}
=a^{(0)}_{jk}u_{\mu_{jk}+\xi_{jk}+1},\\
&&\psi_{a^{(0)},b^{(\xi)}}(z_{jk})u_{\mu_{jk}+\xi_{jk}}
=b^{(\xi)}_{jk}\vep^{\mu_{jk}+\xi_{jk}}u_{\mu_{jk}+\xi_{jk}}=
b^{(0)}_{jk}\vep^{\mu_{jk}}u_{\mu_{jk}+\xi_{jk}}.\qq\qed
\eeqnn

By  Proposition \ref{shift}, we have the following result. 
\begin{pro4}
We consider the representation $(\Phi_{r^{(0)}, s^{(0)},a^{(0)},b^{(\xi)}},\VV)$.
A vector $v\in\VV$ satisfies the condition  
\beq
e_i v=0 \q {\rm for\,\, any\,\,} i=1,\cd, n,
\label{highest-xi}
\eeq
 if and only if $v=cu_{\xi}$ $(c\in\CC)$
where $u_{\xi}\in \VV$ is the fixed basis vector as above. 
\end{pro4}

By this proposition, if we take the parameters properly,
any basis vector $u_{\xi}$ in $\VV$ can be a primitive vector.

\renewcommand{\thesection}{\arabic{section}}
\section{Irreducible $\ufin$-module}
\setcounter{equation}{0}
\renewcommand{\theequation}{\thesection.\arabic{equation}}

Suppose that  parameters $r,s,a,b$ satisfy the conditions (\ref{para1})--(\ref{para3}).
In this case $u_{\vec 0}$ is the unique (up to constant)
primitive vector in $\VV$.
As we defined in the last section, set
 $L(\lm):=\ue u_{\vec 0}$ $(\lm=(\lm_1,\lm_2,\cd,\lm_n),
\,\,0\leq \lm_i\leq l-1)$, 
which is a $\ue$-submodule of $\VV$.
In this section we shall see several properties of this module and 
it amounts to an irreducible $\ufin$-module.

\subsection{Root vectors}
First, we see higher root vectors in $\ue$.
There are several definitions for them.
We shall introduce two of them here and 
discuss their relations.

The first one is defined by using (\ref{h-rt+}) and (\ref{h-rt-}).
We also denote them by $e_{\al}$ and $f_{\al}$ ($\al\in\Delta_+$).

\newtheorem{lm5}{Lemma}[section]
\begin{lm5}
The root vectors $e_{\al}$ and $f_{\al}\in\ue$ ($\al\in\Delta_+$) 
defined by (\ref{h-rt+}) and (\ref{h-rt-}) 
 satisfy the relations (\ref{com-e}), (\ref{com-f}),
(\ref{com1-e}), (\ref{com2-e}), (\ref{com1-f}) and (\ref{com2-f})
in $\ue$.
\end{lm5}

{\sl Proof.}
The proof for the $f_{\al}$ case is similar to the $e_{\al}$ case,
thus we shall only see the $e_{\al}$ case.
We shall show (\ref{com-e}), (\ref{com1-e}) and (\ref{com2-e}) simultaneously
by the induction on the height of roots.
Set $\al=\al_j+\al_{j+1}\cd +\al_k,$ $\beta=\al-\al_j$ and 
$\gamma=\beta-\al_{j+1}\in \Delta_+$ ($j< k$).
If the condition $(\al,\al_i)=0$  and $i<g(\al)(=k)$ hold, we know that 
$i\in \{1,2,\cd,j-2, j+1,j+2,\cd,k-1\}$.
If $i\in\{1,2,\cd,j-2, j+2,\cd,k-1\}$, by the hypothesis of the induction we have
$e_ie_{\beta}=e_{\beta}e_i$ and then
$$
e_ie_{\al}=e_i(\vep^{-1}e_je_{\beta}-e_{\beta}e_j)=(\vep^{-1}e_je_{\beta}-e_{\beta}e_j)e_i
=e_{\al}e_i.
$$
If $i=j+1$, by the hypothesis of the induction we have
$\vep e_{j+1}e_{\beta}=e_{\beta}e_{j+1}$, 
$e_{j}e_{\gamma}=e_{\gamma}e_{j}$ and then
\beqnn
\hspace{-20pt}e_{j+1}e_{\al}&\hspace{-20pt}&=e_{j+1}(\vep^{-1}e_je_{\beta}-e_{\beta}e_j)
=e_{j+1}(\vep^{-1}e_j(\vep^{-1}e_{j+1}e_{\gamma}-e_{\gamma}e_{j+1})
-e_{\beta}e_j)\\
&&=\vep^{-2}e_{j+1}e_je_{j+1}e_{\gamma}-\vep^{-1}e_{j+1}e_je_{\gamma}e_{j+1}
-e_{j+1}e_{\beta}e_j\\
&&=\vep^{-1}e_{\beta}e_{j+1}e_j+e_{\gamma}e_{j+1}e_je_{j+1}+\vep^{-1}e_je_{\beta}e_{j+1}
-e_{\beta}e_je_{j+1}-e_{\gamma}e_{j+1}e_je_{j+1}-\vep^{-1}e_{\beta}e_{j+1}e_j\\
&&=\vep^{-1}e_je_{\beta}e_{j+1}-e_{\beta}e_je_{j+1}=e_{\al}e_{j+1}.
\eeqnn
Here we used the formula; 
\beqnn
e_{j+1}e_je_{j+1}e_{\gamma}& =& 
\vep e_{\beta}e_{j+1}e_j+\vep^2e_{\gamma}e_{j+1}e_je_{j+1}+ \vep e_je_{\beta}e_{j+1},\\
e_{j+1}e_je_{\gamma}e_{j+1}& = &
e_{\beta}e_je_{j+1}+e_{\gamma}e_{j+1}e_je_{j+1}.
\eeqnn

Next, if the condition $(\al,\al_i)=-1$  and $i<g(\al)(=k)$ hold, we know that 
$i=j-1$. In this case, we have
\beqnn
e_{j-1}e_{\al+\al_{j-1}}&&=e_{j-1}(\vep^{-1}e_{j-1}e_{\al}-e_{\al}e_{j-1})\\
&&=e_{j-1}(\vep^{-1}e_{j-1}(\vep^{-1}e_je_{\beta}-e_{\beta}e_j)
-(\vep^{-1}e_je_{\beta}-e_{\beta}e_j)e_{j-1})\\
&&=\vep^{-2}([2]e_{j-1}e_je_{j-1}e_{\beta}-e_je_{j-1}^2e_{\beta})
-\vep^{-1}e_{\beta}e_{j-1}^2e_j-e_{j-1}(\vep^{-1}e_je_{\beta}-e_{\beta}e_j)e_{j-1}\\
&&=((-\vep^{-2}e_je_{\beta}+\vep^{-1}e_{\beta}e_j)e_{j-1}
+\vep^{-2}e_{j-1}(\vep^{-1}e_je_{\beta}-e_{\beta}e_j))e_{j-1}\\
&&=\vep^{-1}(\vep^{-1}e_{j-1}e_{\al}-e_{\al}e_{j-1})e_{j-1}
=\vep^{-1}e_{\al+\al_{j-1}}e_{j-1}.
\eeqnn
 
Finally, under the same condition as above,
 using $e_je_{\al}=\vep^{-1}e_{\al}e_j$, $e_{\beta}e_{\al}=\vep e_{\al}e_{\beta}$,
$e_{j-1}e_{\beta}=e_{\beta}e_{j-1}$ and 
$e_je_{\beta}^2-[2]e_{\beta}e_je_{\beta}+e_{\beta}^2e_j
=e_j^2e_{\beta}-[2]e_je_{\beta}e_j+e_{\beta}e_j^2=0$
we have
\beqnn
e_{\al}^2e_{j-1} &&=(\vep^{-1}e_je_{\beta}-e_{\beta}e_j)e_{\al}e_{j-1}
=e_je_{\al}e_{j-1}e_{\beta}-e_{\beta}e_je_{\al}e_{j-1}\\
&&=\vep^{-1}e_j^2e_{j-1}e_{\beta}^2-e_je_{\beta}e_je_{j-1}e_{\beta}
-\vep^{-1}e_{\beta}e_j^2e_{j-1}e_{\beta}+e_{\beta}e_je_{\beta}e_je_{j-1}\\
&& ={{\vep^{-2}}\over{[2]}}e_j^2e_{j-1}e_{\beta}^2-\vep^{-1}e_{\beta}e_j^2e_{j-1}e_{\beta}^2
+{1\over{[2]}}e_{\beta}^2e_{j}^2e_{j-1}\\
&&=e_{j-1}(-\vep^{-2}e_je_{\beta}e_je_{\beta}+\vep^{-1}e_{\beta}e_j^2e_{\beta}
+\vep^{-1}e_je_{\beta}^2e_j-e_{\beta}e_je_{\beta}e_j)\\
&&+[2](\vep^{-2}e_je_{\beta}e_{j-1}e_je_{\beta}-\vep^{-1}e_je_{\beta}e_{j-1}e_{\beta}e_j
-\vep^{-1}e_{\beta}e_je_{j-1}e_je_{\beta}+e_{\beta}e_je_{j-1}e_{\beta}e_j)\\
&&=-e_{j-1}e_{\al}^2+[2]e_{\al}e_{j-1}e_{\al},
\eeqnn
and then we have $e_{\al}e_{\al+\al_{j-1}}=\vep e_{\al+\al_{j-1}}e_{\al}$.\qed

\vskip10pt
Here we introduce the alternative definition of 
 root vectors {\cite{DJMM}.
For roots $\al=\al_{i}+\al_{i+1}+\cd \al_{j}$ and 
$\beta=\al_{j+1}+\al_{j+2}+\cd +\al_k$
($i<j<k$), we define
\beqn
&\ovl e_{\al+\beta}=\ovl e_{\al} \ovl e_{\beta}-\vep \ovl e_{\beta}\ovl e_{\al},\\
&\ovl f_{\al+\beta}=\ovl f_{\al} \ovl f_{\beta}-\vep^{-1}\ovl f_{\beta}\ovl f_{\al},
\eeqn
where we set $\ovl e_{\al_i}:=e_i$ and $\ovl f_{\al_i}:=f_i$.
Note that this definitions are well-defined, that is, 
these do not depend on the choice of $j$.

We obtain the following simple relations between two types of the root vectors:
\begin{lm5}
\label{root-vec}
For any $\al\in\Delta_+$, we have
\beq
\ovl e_{\al}=\vep^{{\rm height}(\al)-1}e_{\al},\q
\ovl f_{\al}=\vep^{-{\rm height}(\al)+1}f_{\al},
\label{2-root}
\eeq
\beq
\ovl e_{\al}^l=e_{\al}^l,\q
\ovl f_{\al}^l=f_{\al}^l.
\label{l-power}
\eeq
\end{lm5}

{\sl Proof.} The proof of (\ref{2-root}) is done by using induction on the height of roots
and (\ref{l-power}) is immediate consequence of (\ref{2-root}) since $\vep^l=1$. \qed

\subsection{$\ufin$-module structure on $\VV$}
For $\al\in\Delta_+$ we define the actions of $e_{\al}$ and $f_{\al}$ recursively 
by using the formula (\ref{h-rt+}) and (\ref{h-rt-}) as in the
previous subsection.
\newtheorem{pro5}[lm5]{Proposition}
\begin{pro5}
\label{nil}
For any $\al\in\Delta_+$ and $i=1,\cd, n$,  we have
\beq
e_{\al}^l=f_{\al}^l=0 {\hbox{ and }} t_i^{2l}=1{\hbox{ on }}\VV.
\eeq
\end{pro5}

{\sl Proof.}
Since  
$\Phi_{r,s,a,b}(t_i)=\vep^{\lm_i}Z_{in}^2Z_{i-1n}^{-1}Z_{i+1n}^{-1}$ on $\VV$
and $Z_{ij}^l=1$,
it is trivial that $t_i^{2l}=1$. 
To show the nilpotency of $e_{\al}$ and $f_{\al}$ we see the following result 
in \cite{DJMM}.
\begin{pro5}[\cite{DJMM} Proposition 3.4]
For $\al=\al_i+\al_{i+1}+\cd+\al_j$, the actions of $\ovl e_{\al}^l$ 
and $\ovl f_{\al}^l$ on
$\VV$ are given by
\beqn
\hspace{-30pt}\ovl e_{\al}^l
&& ={1\over{(\vep-\vep^{-1})^l}}
\left(\sum_{k_1\geq i,\cd,k_{j-i+1}\geq j}\sum_{p=1}^{j-i+1}
(-1)^{p-1}\theta(k_1\geq \cd \geq k_p<\cd< k_{j-i+1})\right. \nn\\
&& \hspace{-30pt}\left.
C_{i\,k_1\cd k_{p-1}}(C_{i+p-1\, k_p}-C_{i+p-1\,k_p}^{-1})C_{i+p\,k_{p+1}\cd k_{j-i+1}}
D_{i\, k_1\cd k_{j-i+1}} \right) \cdot {\rm id},\nn\\
\hspace{-30pt}\ovl f_{\al}^l&&={1\over{(\vep-\vep^{-1})^l}}
\left(\sum_{k_1\geq n+1-j,\cd,k_{j-i+1}\geq n+1-i}\sum_{p=1}^{j-i+1}
(-1)^{p-1}\theta(k_1\geq \cd \geq k_p<\cd< k_{j-i+1})\right. \nn\\
&&\hspace{-30pt}
\left. \ovl C_{n+1-j\,k_1\cd k_{p-1}}(\ovl C_{n-j+p\, k_p}-\ovl C_{n-j+p\,k_p}^{-1})
\ovl C_{n-j+p+1\,k_{p+1}\cd k_{j-i+1}}
\ovl D_{n-j+1\, k_1\cd k_{j-i+1}}\right)\cdot{\rm id},\nn
\eeqn
where $\theta(X)=1$ if $X$ is true and $\theta(X)=0$ otherwise, and we set 
\beqn
&& C_{i\,k}:=(\vep^{-1}r_ib_{ik}b_{ik-1}b_{i-1k-1}^{-1}b_{i+1k}^{-1})^l,\\
&& \ovl C_{ik}:=(\vep^{-1}s_ib_{k+1-i\,k}b_{k-i\,k-1}
                  b_{k+1-i\,k-1}^{-1}b_{k-i\,k}^{-1})^l,\\
&& D_{ik}:=\prod_{p=k}^n( a_{ip})^l, \,\ovl D_{ik}:=\prod_{p=k}^n(a_{p+1-i\,p}^{-1})^l,
\eeqn
and $\phi_{i\,k_1\cd k_p}:=\phi_{i\,k_1}\cd \phi_{i\,k_p}$ for 
$\phi=C,\ovl C, D,\ovl D$.
\end{pro5}

Applying the specializations of the parameters (\ref{para1}), (\ref{para2})
and (\ref{para3}) to $C_{ik}$, $\ovl C_{ik}$, $D_{ik}$ and $\ovl D_{ik}$,
we have $C_{ik}=\ovl C_{ik}=1$, which implies that 
$C_{n-j+p\, k_p}-C_{n-j+p\,k_p}^{-1}=\ovl C_{n-j+p\, k_p}-\ovl C_{n-j+p\,k_p}^{-1}=0$ 
and then
$\ovl e_{\al}^l=\ovl f_{\al}^l=0$.
Since we have $e_{\al}^l=\ovl e_{\al}^l$ and $f_{\al}^l=\ovl f_{\al}^l$ 
by Lemma \ref{root-vec}, we obtain that 
$e_{\al}^l=f_{\al}^l=0$ on $\VV$.   \qed

\newtheorem{thm5}[lm5]{Theorem}
\begin{thm5}
\label{main}
\begin{enumerate}
\item
If we define the actions of $e_{\al}$ and $f_{\al}$ $(\al\in\Delta_+)$ 
by using (\ref{h-rt+}) and (\ref{h-rt-}), the vector space $\VV$ becomes $\ufin$-module.
\item
The subspace $L(\lm)$
$(\lm=(\lm_1,\cd,\lm_n),\, 
\lm_i\in\{0,1,\cd,l-1\})$ is the irreducible $\ufin$-submodule of
$\VV$.
\end{enumerate}
\end{thm5}

{\sl Proof.}
To show the former half of the theorem, it suffices to check the relations
(\ref{rel1})--(\ref{11}) in Proposition \ref{ufin}.
The relations (\ref{rel1})--(\ref{rel4}) are satisfied since 
$\VV$ is originally $\ue$-module.
The relations (\ref{com-e})--(\ref{com2-f}) are obtained from 
Lemma \ref{root-vec}.
We have the relations (\ref{00}) and (\ref{11}) from Proposition
\ref{nil}.
Thus, we have the well-defined actions of $\ufin$ on $\VV$.

\subsection{Proof of irreducibility}
In orderto show the irreducibility of $L(\lm)$,
we need the following:

\begin{pro5}
\label{pri}
Any finite dimensional $\ufin$-module contains a primitive vector.
\end{pro5}

To show the proposition, we shall show the following lemma
\begin{lm5}
\label{eee0}
Let $L>0$ be a sufficiently large integer. For any $i_1,i_2,\cd,i_L\in I$
we have in $\ufin$,
\beq
e_{i_L}\cd e_{i_2}e_{i_1}=0.
\eeq
\end{lm5}

{\sl Proof.\,\,}
We define a $\ZZ$-gradation on $(\ufin)^+$ by the following way:
As we have seen in Proposition \ref{PBW}, $(\ufin)^+$ has the basis
$$
\{e_{\beta_{N}}^{r_N}e_{\beta_{N-1}}^{r_{N-1}}\cd 
e_{\beta_{1}}^{r_1}\}_{0\leq r_1,\cd,r_N<l}.
$$
Using this, we define
\beq
(\ufin)^+_{d}:=\bigoplus_{r_1{\rm ht}(\beta_1)+\cd+r_N{\rm ht}(\beta_N)=d}\CC
e_{\beta_{N}}^{r_N}e_{\beta_{N-1}}^{r_{N-1}}\cd 
e_{\beta_{1}}^{r_1},
\eeq
where ${\rm ht}(\beta)$ is the height of a root $\beta\in \Delta_+$.
We have 
$$
(\ufin)^+=\bigoplus_d(\ufin)^+_d.
$$
An element in $(\ufin)^+_d$ is called a {\it homogeneous element of degree }$d$.
Since all the relations in $(\ufin)^+$, that is, 
(\ref{com-e}), (\ref{h-rt+}), (\ref{com1-e}, (\ref{com2-e}) and (\ref{00})
are homogeneous, it is well-defined and then we obtain 
$(\ufin)^+_d(\ufin)^+_e\subset (\ufin)^+_{d+e}$ for $d,e\in\ZZ_{\geq0}$.
Hence, $e_{i_L}\cd e_{i_2}e_{i_1}$ is a homogeneous
element of degree $L$.
It immediately follows from Proposition \ref{PBW} that 
the maximum degree is $(l-1)\sum_{i=1}^N{\rm ht}(\beta_i):=J$,
which implies that if $L>J$, $(\ufin)^+_L=0$.
Thus, if $L$ is sufficiently large, a homogeneous element 
$e_{i_L}\cd e_{i_2}e_{i_1}$ must vanish.\qed

\vskip10pt
{\sl Proof of Proposition \ref{pri}.\,\,}
Suppose that a finite dimensional $\ufin$-module
$V$ does not have any primitive vector.
So any non-zero $v\in V$ there exists an infinite sequence 
$i_1,i_2,\cd, i_k,\cd$ $(i_j\in I)$ such that 
all vectors $v$, $e_{i_1}v$, $e_{i_2}e_{i_1}v,
\cd, e_{i_k}\cd e_{i_2}e_{i_1}v,\cd$ never vanish.
But this contradicts to Lemma \ref{eee0}.
Therefore, $V$ contains a primitive vector.\qed

\vskip10pt

Let $W$ be a non-zero submodule of $L(\lm)$. 
By Proposition \ref{pri}, 
$W$ contains a primitive vector.
By the uniqueness of the primitive vector in $\VV$ 
(Proposition \ref{highest-vec}), 
$W$ has to contain $u_{\vec 0}$.
Therefore, $W=L(\lm)$ and then $L(\lm)$ is irreducible.
Here we completed the proof of Theorem \ref{main}(ii).
\qed

\vskip10pt
By Proposition 7.1 in \cite{L2}, for a '$l$-restricted weight' (see 3.2)
$\lm\in P_+$ 
the $\urese$-module $V_{\vep}^{\rm res}(\lm)$ (see 3.2) is identified 
with the irreducible $\ufin$-module, which is isomorphic to $L(\lm)$.
Accodingly, by Theorem \ref{main} we realize the irreducible 
highest weight $\urese$-module $V_{\vep}^{\rm res}(\lm)$
with the $l$-restricted highest weight $\lm$ in the vector space $\VV$.

Our further problem is to write down a basis of $L(\lm)$ explicitly.
Any basis vector of the irreducible $\uq$-module $V(\lm)$ is 
parametrized by "Young tableaux" of shape $\lm$.
So by the constrcution of $V^{\rm res}_{\vep}(\lm)$ in
3.2, we deduce that a basis vector of $L(\lm)$ would be parametrized by
'restricted'(in some sense) Young tableaux.
We would also like to see the structure of the quotient module $\VV/L(\lm)$ or 
the tensor product of $\VV\ot \VV$.

In \cite{Sch}, for the $B_n$, $C_n$ and $D_n$-cases
the analogous presentations of the maximal cyclic representations are given explicitly.
Thus, we can apply the procedure adopted here to them and might hope to obtain the 
irreducible $\ufin(B_n)$ ($\ufin(C_n)$, $\ufin(D_n)$)-modules, 
which will be discussed elsewhere.

\end{document}